\documentclass[12pt, a4paper]{article}

	\usepackage{amsmath, amsthm, amssymb, amsfonts}
	\usepackage[rm,bf,compact,topmarks,calcwidth,pagestyles]{titlesec}
	\usepackage{titletoc}
	\usepackage{multicol}
	\usepackage{bbm}
	\usepackage{setspace}
	\usepackage{graphicx}
	\usepackage[bf,small,center]{caption}

	\newtheorem{thm}{Theorem}
	\newtheorem{prop}[thm]{Proposition}
	\newtheorem{cor}{Corollary}

	\theoremstyle{remark}
	\newtheorem{rem}{Remark}
	
	\newcommand{\gs}{\sigma}

	\newcommand{\gz}{\zeta}

	\newcommand{\Exp}{\mathbb{E}}
	
	\newcommand{\ito}{It\^o}

	\newcommand{\ind}{\mathbbm{1}}      

	\renewcommand{\Pr}{\mathbb{P}}

	\newcommand{\partf}[2]{\frac{\partial{#1}}{\partial{#2}}}
	\newcommand{\partfd}[3]{\frac{\partial^{#3} #1}{\partial #2^{#3} }}

	\newcommand{\Qr}{\mathbb{Q}}

	\newcommand{\Eh}{\hat{\Exp}_x}

	\renewcommand{\cal}[1]{{\mathcal #1}}

\begin{document}

	\title{Bounds for the Transition Density of Time-Homogeneous Diffusion Processes}
	\author{Andrew N. Downes\footnote{Department of Mathematics and Statistics, University of Melbourne, a.downes@ms.unimelb.edu.au}}
	\date{}

	\maketitle
	
	\begin{abstract}
		The paper presents new simple sharp bounds for transition density functions for time-homogeneous diffusions processes. The bounds are obtained under mild conditions on the drift and diffusion coefficients, extending and substantially improving previous results in the literature which were limited to drifts satisfying a linear growth condition. They lead to an asymptotic expression for the time $t$ transition density as $t \rightarrow 0$. While the focus is on the one-dimensional case, an extension to multiple dimensions is discussed. Results are illustrated by numerical examples.
	\end{abstract}
	
\section{Introduction}

	Let $(U_t)$ be a one-dimensional time-homogeneous diffusion process satisfying the stochastic differential equation
		\[	dU_t = \nu(U_t) dt + \gs(U_t) dW_t, \qquad U_0 = u_0,	\]
	where $(W_t)$ denotes a Brownian motion.	The aim of this note is to bound, from above and below, the transition probability density function for $(U_t)$, $p_U(t, u_0, w) := \frac{d}{dw} P_{U}(t,u_0,w)$, where
$P_{U}(t,u_0,w) := \Pr(U_t \leq w | U_0 = u_0)$. While the focus is on the one-dimensional case, the results are easily extended to some special cases in $\mathbb{R}^n$, $n \geq 2$ (see remark at the end of this section). Some simple bounds for the distribution function are also considered.

	Except for a few special cases, the transition functions are unknown for general diffusion processes, so finding approximations to them is an important alternative approach. We use Girsanov's theorem and then a transformation of the Radon-Nikodym density of the type suggested in \cite{Baldi_etal_0802} to relate probabilities for a general diffusion $(U_t)$ to those of a `reference diffusion'. Using a reference diffusion with known transition functions, we are able to derive various bounds for the transition functions under mild conditions on the original process. The results have a simple form and are readily evaluated.
	
	As an aside, the generator of the diffusion $(U_t)$ is given by
		\[	Af(x) = \nu(x) \partf{f}{x} + \frac{1}{2} \gs^2(x) \partfd{f}{x}{2},	\]
	and the transition probability density function is the minimal fundamental solution to the parabolic equation
		\[	\left(A - \partf{}{t}\right)u(t,x) = 0.	\]
	Thus the results presented here also bound solutions to certain types of parabolic partial differential equations.

	Several papers on this topic are available in the literature, especially for bounding the transition probability density. Most recently, \cite{Qian_etal_0504} proposed upper and lower bounds for diffusions whose drift satisfied a linear growth constraint. This appears to be the first such paper to relax the assumption of a bounded drift term. The results in \cite{Qian_etal_0504} will be compared with those obtained in the current paper, although the former can not be used for processes not satisfying the linear growth constraint. To the best of our knowledge, the bounds presented in the current paper are the only ones to relax this constraint, and also appear to generally offer a tightening of the bounds previously available. For further background on diffusions with bounded drift, see e.g.\ \cite{Qian_etal_1103} and references therein.

	In addition, the same ideas allow us to obtain bounds for other functions related to the diffusions. This is not the focus of this note and is not discussed in great detail here, but as an example at the end of Section~\ref{sec:main_result} we consider the density of the process and its first crossing time. This has application in many areas, such as the pricing of financial barrier options. Bounds for other probabilities may be derived in the same manner.

	Consider a one-dimensional time-homogeneous non-explosive diffusion $(U_t)$ governed by the stochastic differential equation (SDE)
		\begin{align}
		\label{eq:original_diff}
			dU_t = \nu (U_t) dt + \gs(U_t) dW_t,
		\end{align}
	where $(W_t)$ is a Brownian motion and $\gs(y)$ is differentiable and non-zero inside the diffusion interval (that is, the the smallest interval $I \subseteq \mathbb{R}$ such that $U_t \in I$ a.s.). As is well-known, one can transform the process to one with unit diffusion coefficient by letting
		\begin{align}
		\label{eq:fn_transform}
			F(y) := \int_{y_0}^y \frac{1}{\gs(u)} du
		\end{align}
	for some $y_0$ from the diffusion interval of $(U_t)$ and then considering $X_t := F(U_t)$ (see e.g.\ \cite{Rogers_xx85}, p.161). By \ito's formula, $(X_t)$ will have unit diffusion coefficient and a drift coefficient $\mu(y)$ given by the composition
		\[	\mu(y) := \left( \frac{\nu}{\gs} - \frac{1}{2} \gs'\right) \circ F^{-1}(y).	\]
	From here on we work with the transformed diffusion process $(X_t)$ governed by the SDE
		\begin{align*}
			dX_t = \mu(X_t) dt + dW_t, \qquad X_0 = F(U_0) =: x.
		\end{align*}
	Conditions mentioned throughout refer to the transformed process $(X_t)$ and its drift coefficient $\mu$.
	
	We will consider the following two cases only (the results extend to diffusions with other diffusion intervals with one finite endpoint by employing appropriate transforms):
		\begin{enumerate}
			\item[] [A] The diffusion interval of $(X_t)$ is the whole real line $\mathbb{R}$.
			\item[] [B] The diffusion interval of $(X_t)$ is $(0, \infty)$.
		\end{enumerate}
	The results extend to diffusions with other diffusion intervals with one finite endpoint by employing appropriate transforms.
	
	For the diffusion $(X_t)$ we will need a reference diffusion $(Y_t)$ with certain characteristics. The reference diffusion must have the same diffusion interval as $(X_t)$ and a unit diffusion coefficient, so that Girsanov's theorem may be applied to $(X_t)$. To be of any practical use, the reference process must also have known transition functions. In case [A], we use the Brownian motion as the reference process, while in case [B] we use the Bessel process of an arbitrary dimension $d \geq 3$.
	
	Recall the definition of the Bessel process $(R_t)$ of dimension $d = 3, 4, \ldots,$ starting at a point $x>0$. This process gives the Euclidean norm of the $d$-dimensional Brownian motion originating at $(x,0, \ldots, 0)$, that is,
		\[  R_t = \sqrt{\bigl(x +W_t^{(1)}\bigr)^2 + \cdots + \bigl(W_t^{(d)}\bigr)^2},	\]
	where the $\bigl(W_t^{(i)}\bigr)$ are independent standard Brownian motions, $i = 1, \ldots, d$. As is well known (see e.g.\ \cite{Revuz_etal_xx99}, p.445), $(R_t)$ satisfies the SDE
		\begin{align*}
			dR_t = \frac{d-1}{2} \frac{1}{R_t}dt + dW_t.
		\end{align*}
	Note that for non-integer values of $d$ the Bessel process of `dimension' $d$ is defined using the above SDE. The process has the transition density function
		\begin{align*}
			p_R(t,y,z) = z \left(\frac{z}{y}\right)^{\eta} t^{-1} e^{-(y^2 + z^2)/2t} {\cal I}_{\eta} \left(\frac{yz}{t} \right),
		\end{align*}
	where $\eta = d/2 -1$ and ${\cal I}_{\eta}(z)$ is the modified Bessel function of the first kind. For further information, see Chapter XI in \cite{Revuz_etal_xx99}.
	
	We denote by $\Pr_x$ and $\Exp_x$ probabilities and expectations conditional on the process in question ($(X_t)$ or some other process, which will be obvious from the context) starting at $x$. We work with the natural filtration $\cal{F}_s := \gs\left(X_u: u \leq s\right)$.
	
	Finally, note that the present work can be easily extended to a class of $n$-dimensional diffusions for $n \geq 2$. If $(X_t)$ is an $n$-dimensional diffusion satisfying the SDE
		\[	dX_t = \mu(X_t) dt + dW_t,	\]
	$(W_t)$ being an $n$-dimensional Brownian motion, then the majority of results can be extended assuming $\mu(\cdot)$ is curl-free. The extension is straight-forward, and in this note we shall only concern ourselves with the one-dimensional case.

\section{Main Results}
\label{sec:main_result}

	This section states and proves a result relating probabilities for the diffusion $(X_t)$ to expectations under an appropriate reference measure. In the case [A], the result may be known, and we state it here for completeness. The extension to case [B] is straight-forward. We then apply this proposition to obtain bounds for transition densities and distributions.

	\textbf{Relation to the Reference Process}
	
	We define the functions $G(y)$ and $N(t)$ as follows, according to the diffusion interval of $(X_t)$:
 		\begin{enumerate}
			\item[][A] If the diffusion interval of $(X_t)$ is $\mathbb{R}$, then we define, for some fixed $y_0 \in \mathbb{R}$,
          \begin{equation}
          	\begin{array}{rl}
            	G(y) \!\!\!&:= \displaystyle\int_{y_0}^{y} \mu(z) dz,\\
            	\vphantom{.}\\
							N(t) \!\!\!&:= \displaystyle\int_0^t \left(\mu'(X_u) + \mu^2(X_u)\right) du.
            \end{array}
          \label{eq:bm_G}
        	\end{equation}
			\item[][B] If the diffusion interval of $(X_t)$ is $(0, \infty)$, then we define, for some fixed $d \geq 3$ (the dimension of the reference Bessel process) and $y_0 > 0$,
					\begin{align*}
							G(y) &:= \int_{y_0}^{y} \left(\mu(z) - \frac{d-1}{2z} \right) dz,\\
							N(t) &:= \int_0^t \left( \mu'(X_u) - \frac{(d-1)(d-3)}{4X_u^2} + \mu^2(X_u)  \right) du.
					\end{align*}
			\end{enumerate}
		\begin{rem}
			For diffusions on $(0, \infty)$, the choice of $d$ is arbitrary subject to $d \geq 3$. Therefore this choice can be used to optimise any bounds presented in the next subsection.
		\end{rem}

\begin{prop}
\label{th:gen_int}
	Assume the the drift coefficient $\mu$ of $(X_t)$ is absolutely continuous. Then, for any $A \in \mathcal{F}_t$,
		\[	\Pr_x(A) = \Eh\left[ e^{G(X_t) - G(x)} e^{-(1/2) N(t)}\ind_A \right],	\]
	where $\Eh$ denotes expectation with respect to the law of the reference process.
\end{prop}
\begin{rem}
	In terms of the original process $(U_t)$ defined in \eqref{eq:original_diff}, the condition of absolute continuity of $\mu(y)$ requires $\nu(z)$ and $\gs'(z)$ to be absolutely continuous.
\end{rem}

\begin{proof}

	The proof is a straight-forward application of Girsanov's theorem and its idea is similar to the one used in \cite{Baldi_etal_0802}. We present the proof for case [A], the proof for case [B] is completed similarly (see \cite{Downes_etal_xx08} for the general approach).
	
	Define $\Qr_x$ to be the reference measure such that under $\Qr_x$, $X_0 = x$ and
		\[	dX_s = d\widetilde{W}_s,	\]
	for a $\Qr_x$ Brownian motion $(\widetilde{W}_s)$. Set
		\begin{align*}
			\gz_s &:= \frac{d\Pr_x}{d\Qr_x} = \exp \left\{ \int_0^s \mu (X_u) d\widetilde{W}_u - \frac{1}{2} \int_0^s \mu^2(X_u) du \right\},
		\end{align*}
	so by Girsanov's theorem under $\Pr_x$ we regain the original process $(X_s)$ satisfying
		\[	dX_s = \mu(X_s) ds + dW_s,	\]
	for a $\Pr_x$ Brownian motion $(W_s)$. The regularity conditions allowing this application of Girsanov's theorem are satisfied (see e.g.\ Theorem~7.19 in \cite{Liptser_etal_xx01}), since under both $\Pr_x$ and $\Qr_x$ the process $(X_s)$ is non-explosive and $\mu(y)$ is locally bounded so we have, for any $t>0$,
	  \[	\Pr_x \left( \int_0^t \mu^2(X_s) ds < \infty \right) = \Qr_x \left( \int_0^t \mu^2(X_s)  ds < \infty \right) = 1.	\]
	We then have, under $\Qr_x$, using \ito's formula and \eqref{eq:bm_G},
		\begin{align}
		\label{eq:dG}
			d G (X_s) &= \mu(X_s) dX_s + \frac{1}{2} \mu'(X_s) (dX_s)^2\notag\\
				&= \mu(X_s) d\widetilde{W}_s + \frac{1}{2} \mu'(X_s) ds.
		\end{align}
	Note that in order to apply \ito's formula, we only require $\mu$ to be absolutely continuous with Radon-Nikodym derivative $\mu'$ (see e.g.\ Theorem~19.5 in \cite{Kallenberg_xx97}). This also implies the above is defined uniquely only up to a set of Lebesgue measure zero, and we are free to assign an arbitrary value to $\mu'$ at points of discontinuity.
	
	Rearranging \eqref{eq:dG} gives
		\[	\int_0^s \mu (X_u) d\widetilde{W}_u = G(X_s) - G(X_0) - \frac{1}{2} \int_0^s \mu'(X_u) du.	\]
	Hence
		\[	\gz_s = \exp \left\{ G(X_s) - G(X_0) - \frac{1}{2} \int_0^s \left( \mu'(X_u) + \mu^2 (X_u) \right) du \right\},	\]
	which together with
		\begin{align*}
			\Pr_x(A) = \Exp_x[\ind_A] = \int \ind_A d\Pr_x = \int \ind_A \gz_t d\Qr_x = \Eh [ \gz_t \ind_A ],
		\end{align*}
	completes the proof of the proposition.
	
\end{proof}

	\textbf{Bounds for Transition Densities and Distributions}

	Define $L$ and $M$ as follows, according to the diffusion interval of $(X_t)$:
 		\begin{enumerate}
			\item[][A] If the diffusion interval of $(X_t)$ is $\mathbb{R}$, then
					\begin{align*}
						L &:= \displaystyle \mbox{ess sup}\left(\mu'(y) + \mu^2(y)\right),\\
						M &:= \displaystyle \mbox{ess inf}\left(\mu'(y) + \mu^2(y)\right),
					\end{align*}
				where the essential supremum/infimum is taken over $\mathbb{R}$.
			\item[][B] If the diffusion interval of $(X_t)$ is $(0, \infty)$, then, for some fixed $d \geq 3$ (the dimension of the reference Bessel process), we put
					\begin{align*}
						L &:= \displaystyle \mbox{ess sup} \left( \mu'(y) - \frac{(d-1)(d-3)}{4y^2} + \mu^2(y)  \right),\\
						M &:= \displaystyle \mbox{ess inf} \left( \mu'(y) - \frac{(d-1)(d-3)}{4y^2} + \mu^2(y)  \right),
					\end{align*}
				where the essential supremum/infimum is taken over $(0, \infty)$.
		\end{enumerate}

	Note that in what follows, in the case [B], the dimension of the reference Bessel process may be chosen so  as to optimise the particular bound. Recall also that $(Y_t)$ denotes the reference process (the Weiner process in case [A], the $d$-dimensional Bessel process in case [B]).

	\begin{cor}
	\label{cor:trans_dens}
			The transition density of the diffusion $(X_t)$ is bounded according to
				\begin{align}
				\label{eq:trans_dens}
					e^{-tL/2} \leq \frac{p_X(t,x,w)}{e^{G(w) - G(x)} p_Y(t,x,w)} \leq e^{-tM/2}.
				\end{align}
	\end{cor}
	\begin{rem}
		The bound is sharp: for a constant drift coefficient $\mu$, equalities hold in \eqref{eq:trans_dens}.
	\end{rem}	
	\begin{proof}
		Recall (see the proof of Proposition~\ref{th:gen_int}) we only required $\mu$ to be absolutely continuous, and its value on a set of Lebesgue measure zero is irrelevant. Hence $L$ (respectively $M$) gives an upper (lower) bound for the integrand in $N(t)$ for all paths. Applying Proposition~\ref{th:gen_int} with $A = \{X_t \in [w, w+h)\}$, $h>0$, gives
			\begin{align*}
				\inf_{w \leq y \leq w+h} e^{G(y) - G(x)} e^{-tL/2} \Pr_x(Y_t \in [w, &w+h)) \leq \Pr_x(X_t \in [w, w+h))\\
				&\leq \sup_{w \leq y \leq w+h} e^{G(y) - G(x)} e^{-tM/2} \Pr_x(Y_t \in [w, w+h)).
			\end{align*}
		Taking the limits as $h \rightarrow 0$ gives the required result.
	\end{proof}
	In the case of bounded $L$ and $M$ this immediately gives an asymptotic expression for the density $p_X(t,x,w)$ as $t \rightarrow 0$.
	\begin{cor}
		If $-\infty < L, M < \infty$, then, as $t \rightarrow 0$,
			\[	p_X(t,x,w) \sim e^{G(w) - G(x)} p_Y(t,x,w),	\]
		uniformly in $x$, $w$.
	\end{cor}

	While the tightest bounds for the transition distribution are obtained by integrating the bounds for the density given above, this does not, in general, yield a simple closed form expression. We mention other, less tight bounds that are simple and are obtained by a further application of Proposition~\ref{th:gen_int}.

	\begin{cor}
	\label{cor:trans_dist}
		The transition distribution function of the diffusion $(X_t)$ admits the following bound: for any $w \in \mathbb{R}$,
			\begin{align*}
				 \inf_{ \ell \leq y \leq w} e^{G(y) - G(x)} e^{-tL/2} P_Y(t,x,w) \leq P_X(t,x,w)
				 \leq \sup_{\ell \leq y \leq w} e^{G(y) - G(x)} e^{-tM/2} P_Y(t,x,w),
			\end{align*}
		where $\ell$ is the lower bound of the diffusion interval.
	\end{cor}
	
	The assertion of the corollary immediately follows from that of Proposition~\ref{th:gen_int} with $A = \{X_t \leq w \}$.
	
	By considering other events (e.g.\ $A = \{ X_t > w \}$), other similar bounds can be derived.
	
	\textbf{Further Probabilities}
	
	While the focus of this note is on bounds for the transition functions, Proposition \ref{th:gen_int} can be used to obtain other useful results. For example, consider
		\[	\eta_X(t, x, y, w)	:= \frac{d}{dw} \Pr_x \left( \sup_{0 \leq s \leq t} X_s \geq y, X_t \leq w \right).	\]
	Such a function has applications in many areas, for example the pricing of barrier options in financial markets. Using ideas similar to the proof of Corollary~\ref{cor:trans_dens} immediately gives
		\begin{cor}
		\label{cor:other_probs}
			For the diffusion $(X_t)$,
			\[	e^{ -tL/2} \leq \frac{\eta_X(t,x,y,w)}{e^{G(w) - G(x)} \eta_Y(t,x,y,w)} \leq e^{ -tM/2}.	\]
		\end{cor}
	Note that for such probabilities the bounds may be improved, if desired, by replacing $L$ and $M$ with appropriate constants on a case-by-case basis. For example, if we are considering the probability our diffusion stays between two constant boundaries at the levels $c_1 < c_2$, then the supremum (for $L$) and infimum (for $M$) need only be taken over the range $c_1 \leq y \leq c_2$.
	
	Other probabilities may be considered in a similar way.
	
\section{Numerical Results}
\label{sec:num_results}

	Here we illustrate the precision of the results from the previous section for transition densities. Bounds from Corollary~\ref{cor:trans_dens} are compared with known transition density functions and previously available bounds for the Ornstein-Uhlenbeck process in the case [A]. For the case [B], we only compare the bounds obtained in the current paper with exact results, since there appears to be no other known bounds in the literature. We also construct a `truncated Ornstein-Uhlenbeck' process in order to compare our results with other bounds available in the literature. For the Ornstein-Uhlenbeck process we also consider an example to illustrate Corollary~\ref{cor:other_probs}.

	\textbf{The Ornstein-Uhlenbeck Process}

\begin{figure}[!htb]
	\begin{centering}
			\includegraphics[width = 0.90 \textwidth, height = 2.5 in]{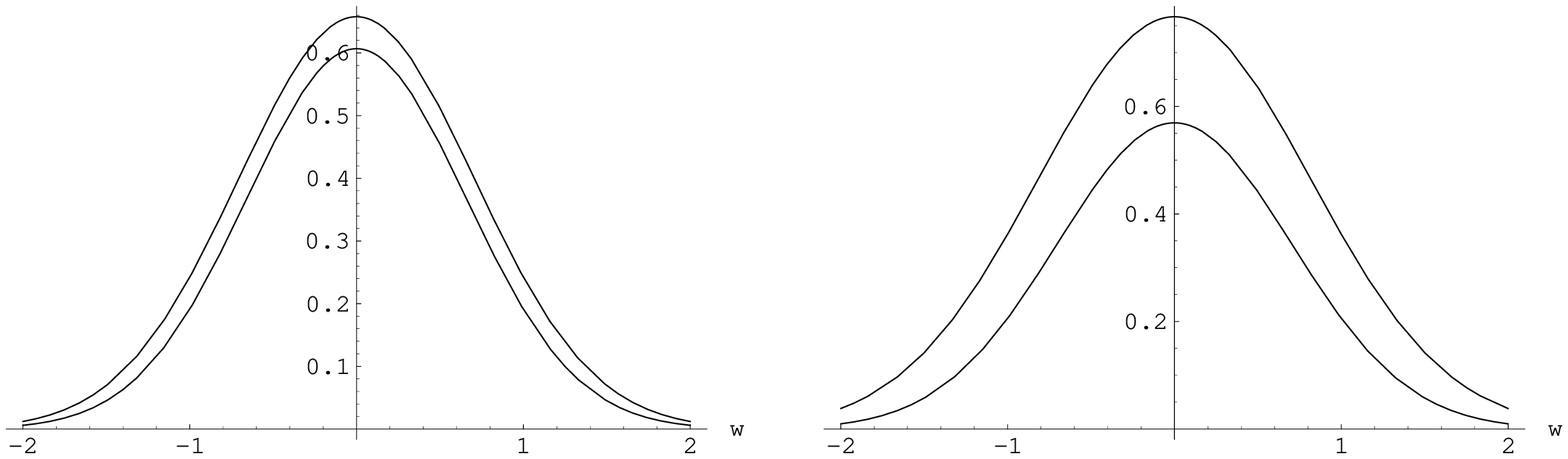}
		\caption{Transition density for an Ornstein-Uhlenbeck process, alongside its upper bound, with $x=0$. The left-hand side displays the functions for $t=1$, the right for $t=2$.}
		\label{fig:OU_dens_cent}
	\end{centering}
\end{figure}

	We consider an Ornstein-Uhlenbeck process $(S_t)$, which satisfies the SDE
		\[	dS_t = -S_t dt + dW_t.	\]
	This process has the transition density
		\[	p_S(t,x,w) = \frac{ e^{t}}{\sqrt{\pi (e^{2 t}-1)}} \exp\left(\frac{\left(w e^{t} - x\right)^2}{1-e^{2 t}} \right),	\]
	see e.g.\ (1.0.6) in \cite{Borodin_etal_xx02}, p.522, and we begin by comparing this with the bound obtained in Corollary~\ref{cor:trans_dens}. Since $\mu(z) = -z$, we have
		\[	M = -1, \qquad G(w) - G(x) = \frac{1}{2}(x^2 - w^2),	\]
	giving the estimate
		\[	p_S(t,x,w) \leq e^{\frac{1}{2}(x^2 - w^2 + t)} p_W(t,x,w).	\]
	Clearly in this case the bound is tighter for smaller values of $|x|$ and $t$. Figure~\ref{fig:OU_dens_cent} displays a plot of the right-hand side of this bound together with the exact density for $x=0$ and $t=1,2$.

\begin{figure}[!htb]
	\begin{centering}
			\includegraphics[width = 0.45 \textwidth, height = 2.5 in]{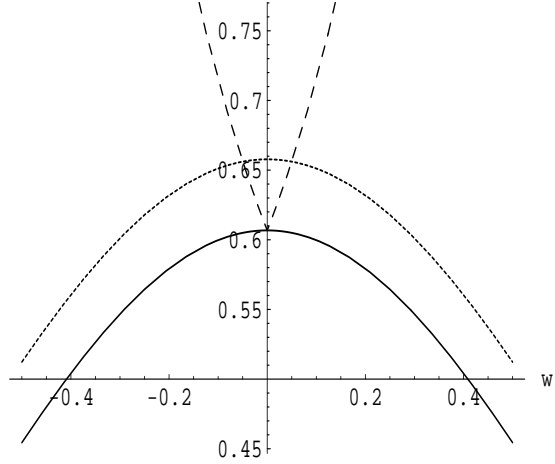}
		\caption{Transition density for an Ornstein-Uhlenbeck process (solid line) compared to bounds given in \cite{Qian_etal_0504} (dashed line) and Corollary~\ref{cor:trans_dens} (dotted line), for $x=0$ and $t=1$.}
		\label{fig:OU_dens_comp}
	\end{centering}
\end{figure}

	To compare our results with other known bounds for transition functions, we look at the bound given by (3.3) in \cite{Qian_etal_0504} (which, to the best of the author's knowledge, is the only bound available for such a process). Figure~\ref{fig:OU_dens_comp} compares this bound with that obtained in Corollary~\ref{cor:trans_dens} and the exact transition density. The values $x=0$ and $t=1$ are used (for the bound in \cite{Qian_etal_0504}, $q=1.2$ seemed to give the best result, see \cite{Qian_etal_0504} for further information on notation). Note that \cite{Qian_etal_0504} gives a sharper bound for $w$ close to zero, but quickly grows to very large values as $|w|$ increases, and in general the bounds presented in this note offer a large improvement. This is typical for all values of $t$, with the effect becoming more pronounced as $t$ decreases. A meaningful lower bound for this process is unavailable by the methods of the present paper, since $L = -\infty$.

\begin{figure}[!htb]
	\begin{centering}
			\includegraphics[width = 0.45 \textwidth, height = 2.5 in]{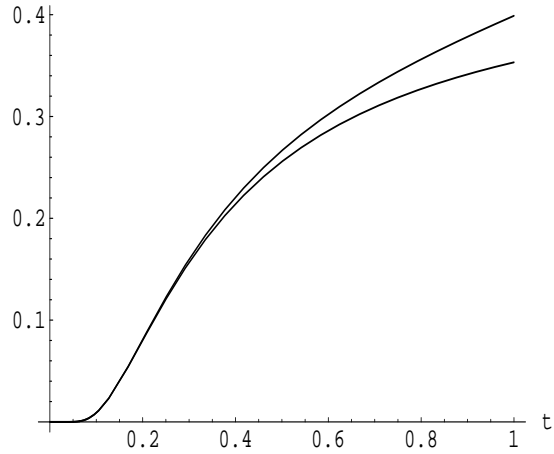}
		\caption{Bound for $\eta_S(t,x,0,z)$ compared with its true value, for $x = z = -0.5$.}
		\label{fig:ou_eta}
	\end{centering}
\end{figure}

	For this example, we briefly look at the bound obtained in Corollary~\ref{cor:other_probs}. We have, see e.g.\ (1.1.8) in \cite{Borodin_etal_xx02}, p. 522, 
		\[	\eta_S(t,x,0,z) = \frac{1}{\sqrt{\pi (1 - e^{-2 t})}} \exp \left( -\frac{(|z| - x e^{-t})^2}{1 - e^{-2 t}} \right).	\]
	Figure~\ref{fig:ou_eta} compares this as a function of $t \in [0,1]$ with the bound obtained in Corollary~\ref{cor:other_probs},
		\begin{align*}
			\eta_S(t,x,0,z) &\leq \exp \left\{\frac{1}{2} (x^2 - z^2 + t) \right\} \eta_W(t,x,0,z)\\
					&= \exp \left\{ \frac{1}{2} (x^2 - z^2 + t) \right\} \frac{1}{\sqrt{2 \pi t}} \exp \left\{ - \frac{1}{2t} (|z| - x)^2 \right\},
		\end{align*}
	where $\eta_W(t,x,0,z)$ is given by (1.1.8) on p. 154 of \cite{Borodin_etal_xx02}.
	
	\textbf{The Truncated Ornstein-Uhlenbeck Process}

\begin{figure}[!htb]
	\begin{centering}
			\includegraphics[width = 0.90 \textwidth, height = 2.5 in]{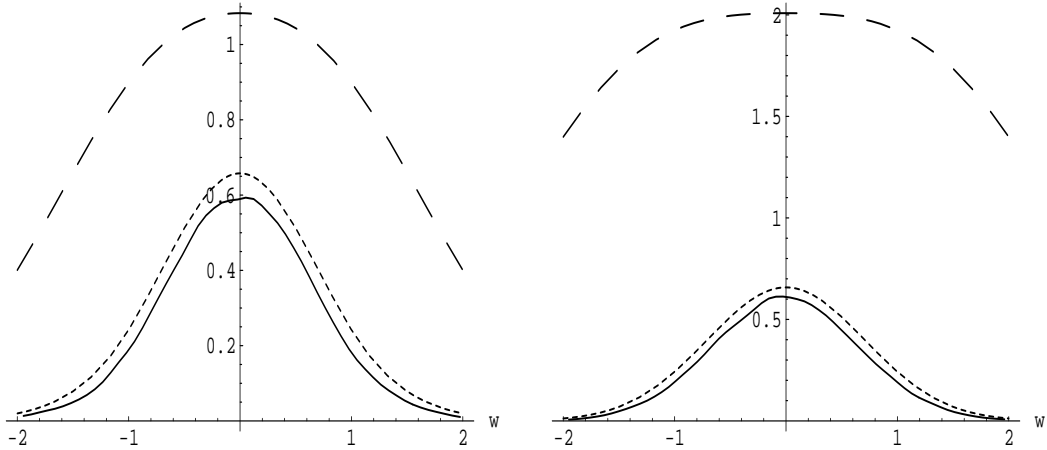}
		\caption{Simulated density and bounds for the transition density of the truncated Ornstein-Uhlenbeck process $(\overline{S}_t)$. The solid lines give the simulated densities, the dotted lines the bound given in Corollary~\ref{cor:trans_dens} and the dashed lines the bounds from \cite{Qian_etal_1103}. Both graphs display the functions for $x=0$ and $t=1$, while the left graph displays them for $c=1$, the right for $c=2$.}
		\label{fig:trunc_ou}
	\end{centering}
\end{figure}

	Other density bounds available in the literature hold only for processes which have bounded drift. For completeness we compare one such bound with the results of this paper. We use the bound in \cite{Qian_etal_1103}, which is the most recent for bounded drift and seems to give the best results over a large domain. To use these results, however, we need a process with bounded drift. As such, we have chosen the `truncated Ornstein-Uhlenbeck' process, which we define as a process $(\overline{S}_t)$ satisfying the SDE
		\[	d\overline{S}_t = \mu(\overline{S}_t) dt + dW_t,	\]
	where, for a fixed $c > 0$,
		\[	\mu(z) = 
					\begin{cases}
						c, & \qquad z < -c,\\
						-z, & \qquad |z| \leq c,\\
						-c, & \qquad z > c.
					\end{cases}
		\]
	For this process we again have $M=-1$ and, assuming $|x| \leq c$,
		\[	G(w) - G(x) =
					\begin{cases}
						\frac{1}{2}(c^2 + x^2)  + cw, & \qquad w < -c,\\
						\frac{1}{2}(x^2 - w^2), & \qquad |w| \leq c,\\
						\frac{1}{2}(c^2 + x^2)  - cw, & \qquad w > c.
					\end{cases}
		\]
	Figure~\ref{fig:trunc_ou} displays the bounds from Corollary~\ref{cor:trans_dens} together with those in \cite{Qian_etal_1103} for different values of $c$ with $x=0$ and $t=1$. Smaller values of $c$ move the bounds closer together, however for the given choice of $x$ and $t$ they do not touch until we use the (rather severe) truncation $c \approx 0.45$. In general the method outlined in this note provides a dramatic improvement. We have also plotted an estimate for the transition density using simulation. The simulation was performed using the predictor-corrector method (see e.g.\ \cite{Kloeden_etal_xx94} p.198), with $10^5$ simulations and $100$ time-steps.
	
	\textbf{A Diffusion on $(0, \infty)$}

	Finally we consider a process from the case [B]. The author believes this is the first paper to present a bound on transition densities without the linear growth constraint. The process $(V_t)$ satisfying the SDE
		\begin{align}
		\label{eq:feller_sde}
			dV_t = (p V_t + q)dt + \sqrt{2 r V_t} dW_t
		\end{align}
	with $p$, $q \in \mathbb{R}$ and $r>0$, has a known transition density (see (26) in \cite{Giorno_etal_0686}). After applying the transform $Z_t = F(V_t)$, with $F(y) = \sqrt{\frac{2}{r}y}$ by \eqref{eq:fn_transform}, we obtain the process
		\[	dZ_t = \mu(Z_t) dt + dW_t,	\]
	where
		\[	\mu(y) = \frac{p}{2} y + \frac{1}{y}\left( \frac{q}{r} - \frac{1}{2} \right).	\]
	For $q > r$ this dominates the drift of a Bessel process of order $2q/r > 2$ so is clearly a diffusion on $(0, \infty)$.
	
	We take the values $q=2.5$, $r=1$ and $p=1$. Using these values, we have
		\[	M = \inf_{0 \leq y \leq \infty} \left[ \frac{y^2}{4} + 2.5 + \frac{1}{y^2} \left(2 - \frac{(d-1)(d-3)}{4}\right) \right],	\]
	and
		\begin{align*}
			G(y) - G(x) = \frac{1}{4}(y^2 - x^2) + c \log \left( \frac{y}{x} \right),
		\end{align*}
	where $d$ is the order of the reference Bessel process and $c = 2 - (d-1)/2$.
		
\begin{figure}[!htb]
	\begin{centering}
			\includegraphics[width = 0.45 \textwidth, height = 2.5 in]{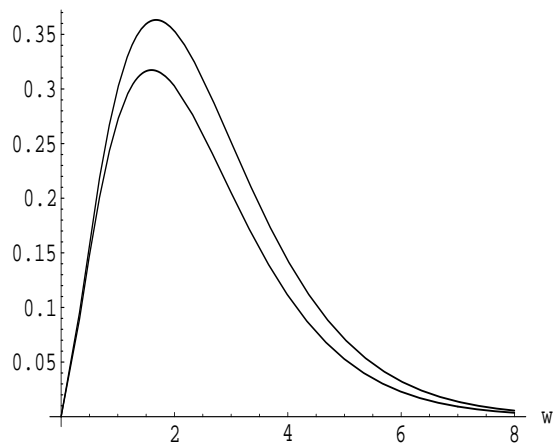}
		\caption{Transition density for the diffusion \eqref{eq:feller_sde}, alongside its upper bound, with $x=0.5$ and $t=0.5$.}
		\label{fig:rplus_dens}
	\end{centering}
	\end{figure}

	It remains to choose the order of the reference Bessel process. It is not clear how to define the `best' order of the reference process for a range of $z$ values, as for fixed $t$ and $x$ the upper bound for $p_Z(t,x,w)$ is minimised for different values of $d$ depending on the value of $z$. In Figure~\ref{fig:rplus_dens} we have taken $t=x=0.5$ and used $d=4.7$, however depending on the relevant criterion improvements can be made. Again, a meaningful lower bound for this process is unavailable by the methods of this paper, since $L= -\infty$.
	
	\newpage

    {\bf Acknowledgements:} This research was supported by the ARC Centre of Excellence for Mathematics and Statistics of Complex Systems. The author is grateful for many useful discussions with K. Borovkov which lead to improvements in the paper.

	\bibliography{BoundaryCrossing,TransDensBounds,GeneralBooks}

\end{document}